\documentclass[aps,prl,preprint]{revtex4}

\begin{document}

\title{When physics helps mathematics: calculation of the sophisticated 
multiple integral}

\author{A.~L.~Kholodenko}
\affiliation{375 H.L.Hunter Laboratories, Clemson University,
Clemson, SC 29634-0973, USA}
\email{string@clemson.edu}

\author{Z.~K.~Silagadze}
\affiliation{Budker Institute of Nuclear Physics 
and Novosibirsk State University, 630 090, Novosibirsk, Russia}
\email{silagadze@inp.nsk.su}

%\date{\today}

\begin{abstract}
There exists a remarkable connection between the quantum mechanical 
Landau-Zener problem and purely classical-mechanical problem of a ball
rolling on a Cornu spiral. This correspondence allows us to calculate
a complicated multiple integral, a kind of multi-dimensional generalization
of Fresnel integrals. A direct method of calculation is also considered but
found to be successful only in some low-dimensional cases. As a byproduct
of this direct method, an interesting new integral representation for
$\zeta(2)$ is obtained.
\end{abstract}

\maketitle

\section{Introduction}
According to Vladimir Arnold \cite{0}, mathematics
can be considered as some branch of physics. In writing this note  we
have no intentions to advocate such point of view. Nevertheless, in our
opinion the following calculus problem is a hard nut to crack if only
mathematical considerations are being used. The problem lies in exactly
calculating  the multiple integral of the type 
\begin{equation}
I_{n}=\int\limits_{-\infty }^{\infty }ds_{1}\int\limits_{-\infty
}^{s_{1}}ds_{2}\cdots \int\limits_{-\infty }^{s_{2n-1}}ds_{2n}\;\cos {%
(s_{1}^{2}-s_{2}^{2})}\;\cdots \cos {(s_{2n-1}^{2}-s_{2n}^{2})}.  \label{eq1}
\end{equation}%
Because of the $s_{1}\leftrightarrow s_{2}$ symmetry, the $n=1$ case is
simple 
\begin{equation}
I_{1}=\frac{1}{2!}\int\limits_{-\infty }^{\infty }ds_{1}\int\limits_{-\infty
}^{\infty }ds_{2}\;\left( \cos {s_{1}^{2}}\,\cos {s_{2}^{2}}+\sin {s_{1}^{2}}%
\,\sin {s_{2}^{2}}\right) =\frac{\pi }{2},  \label{eq2}
\end{equation}%
since its calculation involves known Fresnel integrals 
\[
\int\limits_{-\infty }^{\infty }ds\;\cos {s^{2}}=\int\limits_{-\infty
}^{\infty }ds\;\sin {s^{2}}=\sqrt{\frac{\pi }{2}}.
\]%
However, already for $n\geq 2$ the above symmetry is lost and things quickly
become messy. The $n=2$ and  $n=3$ cases lay at the borderline. They can be 
done with some efforts even though the calculations become noticeably more 
involved. Unfortunately, they do not admit an apparent generalization by using 
the induction method. Already for $n=4$ the attempt to use the same methods 
meets difficulties. Nevertheless, we found a way to calculate such type 
of integrals by invoking some physical arguments.

\section{"If you cannot solve a problem, then there could be an easier
problem you can solve"}

Trying to follow George P\'{o}lya's advice, let us consider the following
system of ordinary differential equations 
\begin{equation}
\frac{d}{ds}\left( 
\begin{array}{c}
x \\ 
y \\ 
z%
\end{array}%
\right) =\frac{1}{R}\left( 
\begin{array}{ccc}
0 & 0 & -\sin {\frac{as^{2}}{2}} \\ 
0 & 0 & \cos {\frac{as^{2}}{2}} \\ 
\sin {\frac{as^{2}}{2}} & -\cos {\frac{as^{2}}{2}} & 0%
\end{array}%
\right) \left( 
\begin{array}{c}
x \\ 
y \\ 
z%
\end{array}%
\right) .  \label{eq3}
\end{equation}%
Such system of equations emerges naturally when one tries to describe a
motion of a sphere $S^{2}$ rolling on the flat surface \textbf{R}$^{2}$
without slippage. More accurately, it describes the rolling of a sphere of
radius $R$ along the Cornu spiral on \textbf{R}$^{2}$ whose curvature $%
\kappa =as$ is proportional to the arc-length $s$ \cite{2,3}. But what is
the relation of this problem to our integral (\ref{eq1})?

As the matrices 
\[
M(s)=\left( 
\begin{array}{ccc}
0 & 0 & -\sin {\frac{as^{2}}{2}} \\ 
0 & 0 & \cos {\frac{as^{2}}{2}} \\ 
\sin {\frac{as^{2}}{2}} & -\cos {\frac{as^{2}}{2}} & 0%
\end{array}%
\right) 
\]%
do not commute for different values of $s$, the solution of (\ref{eq3}) is
given by the time ordered exponential 
\[
U(s,s_{0})=T\exp {\left( \int\limits_{s_{0}}^{s}M(s_{1})\,ds_{1}\right) }%
=1+\int\limits_{s_{0}}^{s}M(s_{1})\,ds_{1}+ 
\]%
\begin{equation}
\;\;\;\;\;\;
\int\limits_{s_{0}}^{s}\,ds_{1}\int\limits_{s_{0}}^{s_{1}}\,ds_{2}\,M(s_{1})%
\,M(s_{2})+\int\limits_{s_{0}}^{s}\,ds_{1}\int\limits_{s_{0}}^{s_{1}}%
\,ds_{2}\int\limits_{s_{0}}^{s_{2}}\,ds_{3}\,M(s_{1})\,M(s_{2})\,M(s_{3})+%
\cdots  \label{eq4}
\end{equation}%
so that 
\begin{equation}
\left( 
\begin{array}{c}
x(s) \\ 
y(s) \\ 
z(s)%
\end{array}%
\right) =U(s,s_{0})\left( 
\begin{array}{c}
x(s_{0}) \\ 
y(s_{0}) \\ 
z(s_{0})%
\end{array}%
\right) .  \label{eq5}
\end{equation}%
Let us take a closer look at the individual terms of the above infinite
series (\ref{eq4}). Writing the matrix $M(s_{i})$ in the block-form 
\[
M_{i}\equiv M(s_{i})=\left( 
\begin{array}{cc}
0 & \chi _{i} \\ 
-\chi _{i}^{T} & 0%
\end{array}%
\right) , 
\]%
where 
\[
\chi _{i}=\frac{1}{R}\left( 
\begin{array}{c}
-\sin {\frac{as_i^2}{2}} \\
\;\;\;\;\cos {\frac{as_i^2}{2}}%
\end{array}%
\right) , 
\]%
it is easy to prove by induction that 
\[
M_{1}M_{2}\cdots M_{2n+1}= 
\]%
\[
\left( 
\begin{array}{cc}
0 & (-1)^{n}\chi _{1}\chi _{2}^{T}\cdots \chi _{2n-1}\chi _{2n}^{T}\chi
_{2n+1} \\ 
(-1)^{n+1}\chi _{1}^{T}\chi _{2}\cdots \chi _{2n-1}^{T}\chi _{2n}\chi
_{2n+1}^{T} & 0%
\end{array}%
\right) , 
\]%
and 
\[
M_{1}M_{2}\cdots M_{2n}=\left( 
\begin{array}{cc}
(-1)^{n}\chi _{1}\chi _{2}^{T}\cdots \chi _{2n-1}\chi _{2n}^{T} & 0 \\ 
0 & (-1)^{n}\chi _{1}^{T}\chi _{2}\cdots \chi _{2n-1}^{T}\chi _{2n}%
\end{array}%
\right) . 
\]%
But 
\[
\chi _{i}^{T}\chi _{i+1}=\frac{1}{R^{2}}\,\cos {\frac{a}{2}%
(s_{i}^{2}-s_{i+1}^{2})}. 
\]%
Besides, then our calculations is supplemented by the following initial and 
boundary conditions 
\begin{equation}
s_{0}=-\infty ,\;\;s=\infty ,\;\;x(-\infty )=y(-\infty )=0,\;\;z(-\infty )=1,
\label{init}
\end{equation}%
in view of these relations and taking into account (\ref{eq5}) and 
(\ref{eq4}), we finally obtain 
\begin{equation}
z(\infty )=1+\sum\limits_{n=1}^{\infty }\frac{(-1)^{n}}{R^{2n}}%
\int\limits_{-\infty }^{\infty }ds_{1}\int\limits_{-\infty
}^{s_{1}}ds_{2}\cdots \int\limits_{-\infty
}^{s_{2n-1}}ds_{2n}\;f_{n}(s_{1},s_{2},\ldots ,s_{2n}),  \label{eq6}
\end{equation}%
where 
\[
f_{n}(s_{1},s_{2},\ldots ,s_{2n})=\cos {\frac{a}{2}(s_{1}^{2}-s_{2}^{2})}%
\,\cos {\frac{a}{2}(s_{3}^{2}-s_{4}^{2})}\cdots \cos {\frac{a}{2}%
(s_{2n-1}^{2}-s_{2n}^{2})}. 
\]%
After rescaling 
\begin{equation}
s_{i}\rightarrow \sqrt{\frac{2}{a}}\,s_{i},  \label{resc}
\end{equation}%
our original integral (\ref{eq1}) indeed shows up in the series (\ref{eq6}): 
\begin{equation}
z(\infty )=1+\sum\limits_{n=1}^{\infty }(-1)^{n}\left( \frac{2}{aR^{2}}%
\right) ^{n}I_{n}.  \label{eq7}
\end{equation}%
Thus, if we can find $z(\infty )$, we can find $I_{n}$ as well! But how can
we find $z(\infty )$?

\section{Spinorization and the Hopf map}

Now we shall follow the advice of Jacques Hadamard: \emph{\textquotedblleft
The shortest path between two islands of truths in the real domain passes
through the complex plane\textquotedblright .} By adopting his suggestion to
our case it is convenient at this point to use two complex variables $a$ and 
$b$ instead of three real variables $x,y$ and $z$ through relations \cite{4} 
\begin{equation}
x=ab^{\ast }+ba^{\ast },\;\;y=i(ab^{\ast }-ba^{\ast }),\;\;z=aa^{\ast
}-bb^{\ast }.  \label{eq8}
\end{equation}%
Notice: introduction of two complex variables is equivalent to looking at
solution for our problem in \textbf{C}$^{2}!$ Furthermore, in view of the
initial conditions and (\ref{eq3}), the variables x,y and z are constrained to
unit sphere $S^{2}.$ This causes the variables $a$ and $b$ to be constrained
to $S^{3}$, that is to obey the equation $\left\vert a\right\vert
^{2}+\left\vert b\right\vert ^{2}=1.$ See Ref.\cite{4} for more details.
Interestingly enough, under these conditions Eq.\ref{eq8} describes the Hopf
map $S^{3}\rightarrow S^{2}$ \cite{5}. From Ref.\cite{4} it can be seen that
the complex variables $a$ and $b$ must satisfy the following system of
differential equations 
\begin{equation}
i\frac{d}{ds}\left( 
\begin{array}{c}
a \\ 
b%
\end{array}%
\right) =-\frac{1}{2R}\left( 
\begin{array}{cc}
0 & e^{-ias^{2}/2} \\ 
e^{ias^{2}/2} & 0%
\end{array}%
\right) \left( 
\begin{array}{c}
a \\ 
b%
\end{array}%
\right) ,  \label{eq9}
\end{equation}%
if the real variables $x,y$ and $z$ to satisfy (\ref{eq3}). Now we can again
formally solve (\ref{eq9}) by using the time-ordered exponential series.
This time, however, the solution is known and it was obtained by Rojo in 
\cite{1}. We just shortly repeat it to ensure the continuity of our
exposition. The solution is formally given by 
\begin{equation}
\left( 
\begin{array}{c}
a(s) \\ 
b(s)%
\end{array}%
\right) =U(s,s_{0})\left( 
\begin{array}{c}
a(s_{0}) \\ 
b(s_{0})%
\end{array}%
\right) ,  \label{eq10}
\end{equation}%
where 
\begin{equation}
U(s,s_{0})=1-i\int\limits_{s_{0}}^{s}H(s_{1})\,ds_{1}+(-i)^{2}\int%
\limits_{s_{0}}^{s}\,ds_{1}\int\limits_{s_{0}}^{s_{1}}\,ds_{2}\,H(s_{1})%
\,H(s_{2})+\cdots  \label{eq11}
\end{equation}%
with 
\[
H(s)=-\frac{1}{2R}\left( 
\begin{array}{cc}
0 & e^{-ias^{2}/2} \\ 
e^{ias^{2}/2} & 0%
\end{array}%
\right) . 
\]%
The imposed initial conditions (\ref{init}) are now translated into 
\begin{equation}
a(-\infty )=1,\;\;\;b(-\infty )=0.  \label{eq12}
\end{equation}%
Since the unitary evolution (\ref{eq10}) conserves the norm $aa^{\ast
}+bb^{\ast }$, from (\ref{eq12}) we obtain back the equation for 3-sphere,
that is $aa^{\ast }+bb^{\ast }=1,$ valid for any \textquotedblleft
time\textquotedblright\ $s$. This result implies that 
\begin{equation}
z(\infty )=|a(\infty )|^{2}-|b(\infty )|^{2}=2|a(\infty )|^{2}-1.
\label{eq13}
\end{equation}%
Evidently, we need only to calculate $a(\infty )$ to obtain $z(\infty )$.

By examining the product of $H(s_{i})$ matrices, we observe that for the odd
number of multipliers the matrix product does not have non-zero diagonal
terms and, hence, does not contribute to $a(\infty )$ thanks to the initial
conditions (\ref{eq12}). The remaining terms with even number of multipliers
have easily calculable non-zero diagonal elements \ so that we get 
\[
a(\infty )=1+\left( \frac{i}{2R}\right) ^{2}\int\limits_{-\infty }^{\infty
}ds_{1}\,e^{-ias_{1}^{2}/2}\int\limits_{-\infty
}^{s_{1}}ds_{2}\,e^{ias_{2}^{2}/2}+\cdots . 
\]%
After rescaling (\ref{resc}) this expression acquires the form 
\begin{equation}
a(\infty )=1+\sum\limits_{n=1}^{\infty }(-1)^{n}\left( \frac{1}{2aR^{2}}%
\right) ^{n}J_{n},  \label{eq14}
\end{equation}%
where 
\begin{equation}
J_{n}=\int\limits_{-\infty }^{\infty
}ds_{1}\,e^{-is_{1}^{2}}\int\limits_{-\infty
}^{s_{1}}ds_{2}\,e^{is_{2}^{2}}\int\limits_{-\infty
}^{s_{2}}ds_{3}\,e^{-is_{3}^{2}}\cdots \int\limits_{-\infty
}^{s_{2n-1}}ds_{2n}\,e^{is_{2n}^{2}}.  \label{eq15}
\end{equation}%
In contrast to $I_{n}$, the multiple integral $J_{n}$ is doable. It can be
calculated as follows \cite{1}. First, we write 
\[
J_{n}=\int\limits_{-\infty }^{\infty }ds_{1}\cdots \int\limits_{-\infty
}^{s_{2n-1}}ds_{2n}\,\theta (s_{1}-s_{2})\cdots \theta
(s_{2n-1}-s_{2n})\,e_{n}(s_{1}^{2},\ldots ,s_{2n}^{2}), 
\]%
with 
\[
e_{n}(s_{1}^{2},\ldots ,s_{2n}^{2})=\exp {%
\{-i(s_{1}^{2}-s_{2}^{2}+s_{3}^{2}-s_{4}^{2}+\cdots
+s_{2n-1}^{2}-s_{2n}^{2})\}}. 
\]%
Then, we use the integral representation for the Heaviside step function 
\begin{equation}
\theta (s)=\frac{1}{2\pi i}\int\limits_{-\infty }^{\infty }d\omega \;\frac{%
e^{i\omega s}}{\omega -i\epsilon }=\left\{ 
\begin{array}{l}
1,\;\;\mathrm{if}\;\;s>0, \\ 
\frac{1}{2},\;\;\mathrm{if}\;\;s=0, \\ 
0,\;\;\mathrm{if}\;\;s<0.%
\end{array}%
\right.  \label{eq16}
\end{equation}%
Using the result, we can perform the integrals over $ds_{i}$. For this
purpose we sequentially complete the squares, e.g. 
\[
i\left[ s_{2}(\omega _{2}-\omega _{1})+s_{2}^{2}\right] =i\left[ \left(
s_{2}+\frac{\omega _{2}-\omega _{1}}{2}\right) ^{2}-\frac{(\omega
_{2}-\omega _{1})^{2}}{4}\right] , 
\]%
and then evaluate the Gaussian integrals 
\[
\int\limits_{-\infty }^{\infty }e^{\pm is^{2}}\,ds=\int\limits_{-\infty
}^{\infty }(\cos {s^{2}}\pm i\sin {s^{2}})\,ds=\sqrt{\frac{\pi }{2}}\,(1\pm
i). 
\]%
As a result, we finally obtain 
\[
J_{n}=\pi ^{n}\int\limits_{-\infty }^{\infty }\frac{d\omega _{1}}{2\pi i}%
\cdots \int\limits_{-\infty }^{\infty }\frac{d\omega _{2n-1}}{2\pi i} 
\]%
\begin{equation}
\frac{\exp {\left\{ \frac{i}{2}\left[ \omega _{2}(\omega _{1}-\omega
_{3})+\omega _{4}(\omega _{3}-\omega _{5})+\cdots +\omega _{2n-2}(\omega
_{2n-3}-\omega _{2n-1})\right] \right\} }}{(\omega _{1}-i\epsilon )(\omega
_{2}-i\epsilon )\cdots (\omega _{2n-1}-i\epsilon )}.  \label{eq17}
\end{equation}%
As it can be seen, the terms quadratic in $\omega _{i}$ are all canceled
thanks to the alternating signs in exponents in (\ref{eq15}). It is this
feature that distinguishes, as we shall demonstrate in the next section, the
calculation of $J_{n}$ from the calculation of $I_{n}$ and makes the
integral $J_{n}$ solvable. For this purpose we rescale even-index variables $%
\omega _{2i}\rightarrow 2\omega _{2i}$ in (\ref{eq17}) and perform integrals
in these variables taking into account (\ref{eq16}). The result is: 
\[
J_{n}=\pi ^{n}\int\limits_{-\infty }^{\infty }\frac{d\omega _{1}}{2\pi i}%
\int\limits_{-\infty }^{\infty }\frac{d\omega _{3}}{2\pi i}\cdots 
\]%
\[
\int\limits_{-\infty }^{\infty }\frac{d\omega _{2n-1}}{2\pi i}\;\frac{\theta
(\omega _{1}-\omega _{3})\theta (\omega _{3}-\omega _{5})\cdots \theta
(\omega _{2n-3}-\omega _{2n-1})}{(\omega _{1}-i\epsilon )(\omega
_{3}-i\epsilon )\cdots (\omega _{2n-1}-i\epsilon )}, 
\]%
or 
\begin{equation}
\;\;J_{n}=\int\limits_{-\infty }^{\infty }\frac{d\omega _{1}}{2\pi i}%
\int\limits_{-\infty }^{\omega _{1}}\frac{d\omega _{3}}{2\pi i}\cdots
\int\limits_{-\infty }^{\omega _{2n-3}}\frac{d\omega _{2n-1}}{2\pi i}\;\frac{%
\pi ^{n}}{(\omega _{1}-i\epsilon )(\omega _{3}-i\epsilon )\cdots (\omega
_{2n-1}-i\epsilon )}.  \label{eq18}
\end{equation}%
It is the symmetry of the integrand in (\ref{eq18}) that makes the
calculation of (\ref{eq18}) as easy as the calculation of $I_{1}$ from
which\ our story had begun. In the present case we obtain 
\begin{equation}
J_{n}=\pi ^{n}\frac{1}{n!}\left[ \int\limits_{-\infty }^{\infty }\frac{%
d\omega _{1}}{2\pi i}\,\frac{1}{\omega _{1}-i\epsilon }\right] ^{n}=\frac{%
\pi ^{n}}{n!}\,[\theta (0)]^{n}=\frac{1}{n!}\left( \frac{\pi }{2}\right)
^{n}.  \label{eq19}
\end{equation}%
Then (\ref{eq14}) shows that 
\begin{equation}
a(\infty )=\exp {\left( -\frac{\pi }{4aR^{2}}\right) },  \label{eq20}
\end{equation}%
and from (\ref{eq13}) we get 
\[
z(\infty )=2\exp {\left( -\frac{\pi }{2aR^{2}}\right) }-1=1+\sum%
\limits_{n=1}^{\infty }(-1)^{n}\frac{2}{n!}\left( \frac{\pi }{2aR^{2}}%
\right) ^{n}. 
\]%
Comparing with (\ref{eq7}), we finally obtain the desired expression for the
integral $I_{n}$: 
\begin{equation}
I_{n}=\frac{2}{n!}\left( \frac{\pi }{4}\right) ^{n}.  \label{eq21}
\end{equation}

\section{Direct calculation of $I_{2}$}

If you are still unhappy by our use of indirect methods of calculation of
deceptively simply looking integral (\ref{eq1}), here we discuss \ some
features of the direct method. Unfortunately, as far as we can see, it works
well only for small values of $n$. \ 

Indeed, let us write 
\begin{equation}
\;\;\;\;\;\;I_{2}=\int\limits_{-\infty }^{\infty }ds_{1}
\int\limits_{-\infty }^{\infty
}ds_{2}\int\limits_{-\infty }^{\infty }ds_{3}\int\limits_{-\infty }^{\infty
}ds_{4}\;\theta _{4}(s_{1},s_{2},s_{3},s_{4})\;\cos {(s_{1}^{2}-s_{2}^{2})}%
\;\cos {(s_{3}^{2}-s_{4}^{2})},  \label{eq22}
\end{equation}%
where 
\begin{equation}
\theta _{4}(s_{1},s_{2},s_{3},s_{4})=\theta (s_{1}-s_{2})\,\theta
(s_{2}-s_{3})\,\theta (s_{3}-s_{4}).  \label{eq23}
\end{equation}%
Then, as in the previous section, we use the integral representation for the
step function and notice that 
\[
\int\limits_{-\infty }^{\infty }ds_{1}\int\limits_{-\infty }^{\infty
}ds_{2}e^{i\omega _{1}s_{1}}e^{-is_{2}(\omega _{1}-\omega _{2})}\;\cos {%
(s_{1}^{2}-s_{2}^{2})}=\frac{1}{2}\int\limits_{-\infty }^{\infty
}ds_{1}\int\limits_{-\infty }^{\infty }ds_{2} 
\]%
\[
\left[ \exp {\left( -i\frac{\omega _{1}^{2}}{4}+i\left( s_{1}+\frac{\omega
_{1}}{2}\right) ^{2}+i\frac{(\omega _{1}-\omega _{2})^{2}}{4}-i\left( s_{2}+%
\frac{\omega _{1}-\omega _{2}}{2}\right) ^{2}\right) }+\right. 
\]%
\[
\left. \exp {\left( i\frac{\omega _{1}^{2}}{4}-i\left( s_{1}-\frac{\omega
_{1}}{2}\right) ^{2}-i\frac{(\omega _{1}-\omega _{2})^{2}}{4}+i\left( s_{2}-%
\frac{\omega _{1}-\omega _{2}}{2}\right) ^{2}\right) }\right] = 
\]%
\begin{equation}
\pi \cos {\left[ \frac{1}{4}\omega _{2}(2\omega _{1}-\omega _{2})\right] }.
\label{eq24}
\end{equation}%
Similar calculations are done for integrals over $ds_{3}$ and $\,ds_{4}$. As
a result, after rescaling $\omega _{2}\rightarrow 2\omega _{2}$, we end up
with the result: 
\begin{equation}
I_{2}=\pi ^{2}\int\limits_{-\infty }^{\infty }\frac{d\omega _{1}}{2\pi i}%
\int\limits_{-\infty }^{\infty }\frac{d\omega _{2}}{2\pi i}%
\int\limits_{-\infty }^{\infty }\frac{d\omega _{3}}{2\pi i}\;\frac{\cos {%
[\omega _{2}(\omega _{1}-\omega _{2})]}\;\cos {[\omega _{2}(\omega
_{3}-\omega _{2})]}}{(\omega _{1}-i\epsilon )(\omega _{2}-i\epsilon )(\omega
_{3}-i\epsilon )}.  \label{eq25}
\end{equation}%
The integrals over $d\omega _{1}$ and $d\omega _{3}$ can be easily
calculated. Indeed, 
\[
\int\limits_{-\infty }^{\infty }\frac{d\omega _{1}}{2\pi i}\,\frac{\cos {%
(\omega _{1}\omega _{2}-\omega _{2}^{2})}}{\omega _{1}-i\epsilon }=\frac{1}{2%
}\int\limits_{-\infty }^{\infty }\frac{d\omega _{1}}{2\pi i}\left[ \frac{%
e^{i\omega _{1}\omega _{2}}}{\omega _{1}-i\epsilon }e^{-i\omega _{2}^{2}}+%
\frac{e^{-i\omega _{1}\omega _{2}}}{\omega _{1}-i\epsilon }e^{i\omega
_{2}^{2}}\right] = 
\]%
\[
\frac{1}{2}\left[ \theta (\omega _{2})e^{-i\omega _{2}^{2}}+\theta (-\omega
_{2})e^{i\omega _{2}^{2}}\right] . 
\]%
Therefore, 
\begin{equation}
I_{2}=\frac{\pi ^{2}}{4}\int\limits_{-\infty }^{\infty }\frac{d\omega _{2}}{%
2\pi i}\,\frac{\left[ \theta (\omega _{2})e^{-i\omega _{2}^{2}}+\theta
(-\omega _{2})e^{i\omega _{2}^{2}}\right] ^{2}}{\omega _{2}-i\epsilon }.
\label{eq26}
\end{equation}%
Now we use the well known result 
\begin{equation}
\frac{1}{\omega -i\epsilon }=P\frac{1}{\omega }+i\pi \delta (\omega )
\label{eq27}
\end{equation}%
to split the integral (\ref{eq26}) into the principal value and the $\delta $%
-function parts: 
\[
I_{2}=\frac{\pi ^{2}}{4}(I_{2P}+I_{2\delta }). 
\]%
Of course, the $\delta $-function part is obtained instantly 
\[
I_{2\delta }=\frac{1}{2}. 
\]%
As for the principal value part, we have 
\[
I_{2P}=\frac{1}{2\pi i}\,\lim_{\epsilon \rightarrow 0}\left[
\int\limits_{\epsilon }^{\infty }\frac{e^{-2i\omega _{2}^{2}}}{\omega _{2}}%
\,d\omega _{2}+\int\limits_{-\infty }^{-\epsilon }\frac{e^{2i\omega _{2}^{2}}%
}{\omega _{2}}\,d\omega _{2}\right] = 
\]%
\[
\hspace*{11mm}\frac{1}{2\pi i}\,\lim_{\epsilon \rightarrow 0}\left[
\int\limits_{\epsilon }^{\infty }\frac{e^{-2i\omega _{2}^{2}}}{\omega _{2}}%
\,d\omega _{2}-\int\limits_{\epsilon }^{\infty }\frac{e^{2i\omega _{2}^{2}}}{%
\omega _{2}}\,d\omega _{2}\right] = 
\]%
\[
-\frac{1}{\pi }\int\limits_{0}^{\infty }\frac{\sin {(2\omega _{2}^{2})}}{%
\omega _{2}}\,d\omega _{2}=-\frac{1}{2\pi }\int\limits_{0}^{\infty }\frac{%
\sin {[(\sqrt{2}\omega _{2})^{2}]}}{(\sqrt{2}\omega _{2})^{2}}\,d[(\sqrt{2}%
\omega _{2})^{2}]=-\frac{1}{4}. 
\]%
In making the last step we have used the Dirichlet integral 
\[
\int\limits_{0}^{\infty }\frac{\sin {\omega }}{\omega }\,d\omega =\frac{\pi 
}{2}. 
\]%
Putting all terms together, we obtain finally 
\begin{equation}
I_{2}=\frac{\pi ^{2}}{16}.  \label{eq28}
\end{equation}%
Which is a special case of (\ref{eq21}), as expected.

Can we think about the general case ($n>2)$ by acting in the manner just
described? By repeating the above steps when $n>2$, we obtain 
\begin{equation}
\;\;\;\;\;\;
I_{n}=\left( \frac{\pi }{2}\right) ^{n}\int\limits_{-\infty }^{\infty }\frac{%
d\omega _{2}}{2\pi i}\int\limits_{-\infty }^{\infty }\frac{d\omega _{4}}{%
2\pi i}\cdots \int\limits_{-\infty }^{\infty }\frac{d\omega _{2n-2}}{2\pi i}%
\;\frac{f_{n}(\omega _{2},\omega _{4},\ldots ,\omega _{2n-2})}{(\omega
_{2}-i\epsilon )(\omega _{4}-i\epsilon )\cdots (\omega _{2n-2}-i\epsilon )},
\label{eq29}
\end{equation}%
where 
\begin{equation}
f_{n}(\omega _{2},\omega _{4},\ldots ,\omega _{2n-2})=\phi (0,\omega
_{2})\phi (\omega _{2},\omega _{4})\cdots \phi (\omega _{2n-4},\omega
_{2n-2})\phi (\omega _{2n-2},0),  \label{eq30}
\end{equation}%
with 
\begin{equation}
\phi (\omega _{1},\omega _{2})=\theta (\omega _{1}-\omega
_{2})\;e^{-i(\omega _{1}^{2}-\omega _{2}^{2})}+\theta (\omega _{2}-\omega
_{1})\;e^{-i(\omega _{2}^{2}-\omega _{1}^{2})}.  \label{eq31}
\end{equation}%
Evidently, for $n>2$ things begin to look rather inconclusive and the above
direct method needs some fresh input in order to be brought to completion.

\section{Direct calculation of $I_3$}

Now let us calculate
$$I_3=\left (\frac{\pi}{2}\right )^3\, I,$$
where
\begin{equation}
I=\int\limits_{-\infty}^\infty \frac{dx}{2\pi i}\,\frac{\phi(0,x)f(x)}
{x-i\epsilon},
\label{eq32}
\end{equation}
and
\begin{equation}
f(x)=\int\limits_{-\infty}^\infty \frac{dy}{2\pi i}\,
\frac{\phi(x,y)\phi(y,0)}{y-i\epsilon}.
\label{eq33}
\end{equation}
Using the relation (\ref{eq27}) in (\ref{eq32}), we get
$$I=\frac{1}{2}\phi(0,0)f(0)+\lim_{\epsilon\to 0}\int\limits_\epsilon^\infty
\frac{dx}{2\pi i}\,\frac{\phi(0,x)f(x)-\phi(0,-x)f(-x)}{x}.$$
However, $\phi(0,0)=1$, while
$$f(0)=\int\limits_{-\infty}^\infty \frac{dy}{2\pi i}\,\frac{\phi(0,y)
\phi(y,0)}{y-i\epsilon}=\frac{1}{2}\phi^2(0,0)+\lim_{\epsilon\to 0}\int
\limits_\epsilon^\infty \frac{dy}{2\pi i}\,\frac{\phi^2(0,y)-\phi^2(0,-y)}
{y},$$
and since $$\phi^2(0,y)-\phi^2(0,-y)=e^{-2iy^2}-e^{2iy^2},$$ when $y\ge
\epsilon>0$, we get
$$f(0)=\frac{1}{2}-\frac{1}{\pi}\int\limits_0^\infty \frac{\sin{(2y^2)}}{y}\,
dy=\frac{1}{4},$$
and, therefore,
\begin{equation}
I=\frac{1}{8}+K,
\label{eq34}
\end{equation}
where
\begin{equation}
K=\lim_{\epsilon\to 0}\int\limits_\epsilon^\infty \frac{dx}{2\pi i}\,
\frac{e^{-ix^2}f(x)-e^{ix^2}f(-x)}{x}.
\label{eq35}
\end{equation}
Now, using again (\ref{eq27}), we have
$$f(x)=\frac{1}{2}\phi(x,0)\phi(0,0)+\lim_{\epsilon\to 0}\int\limits_
\epsilon^\infty \frac{dy}{2\pi i}\,\frac{\phi(x,y)\phi(y,0)-\phi(x,-y)
\phi(-y,0)}{y}$$
and since in (\ref{eq35}) $x>0$, we get
\begin{equation}
\;\;\;\;\;\;f(x)=\frac{1}{2}e^{-ix^2}+\lim_{\epsilon\to 0}\int\limits_
\epsilon^\infty \frac{dy}{2\pi i}\,\frac{\theta(x-y)e^{-ix^2}+
\theta(y-x)e^{-i(2y^2-x^2)}-e^{-i(x^2-2y^2)}}{y}.
\label{eq36}
\end{equation}
Analogously,
\begin{equation}
\;\;\;\;\;\;f(-x)=\frac{1}{2}e^{ix^2}+\lim_{\epsilon\to 0}\int\limits_
\epsilon^\infty \frac{dy}{2\pi i}\,\frac{e^{-i(2y^2-x^2)}-\theta(x-y)
e^{ix^2}-\theta(y-x)e^{-i(x^2-2y^2)}}{y}.
\label{eq37}
\end{equation}
In light of (\ref{eq35}), (\ref{eq36}) and (\ref{eq37}),
$$K=K_1+K_2,$$
where
$$K_1=\frac{1}{2}\int\limits_0^\infty\frac{dx}{2\pi i}\,
\frac{e^{-2ix^2}-e^{2ix^2}}{x}=-\frac{1}{2\pi}\int\limits_0^\infty
\frac{\sin{(2x^2)}}{x}\,dx=-\frac{1}{8}$$
and
$$K_2=-\int\limits_0^\infty\frac{dx}{x}\int\limits_0^\infty
\frac{dy}{y}\,\frac{\theta(x-y)\cos{(2x^2)}+\theta(y-x)\cos{(2y^2)}-
\cos{2(x^2-y^2)}}{2\pi^2}.$$
After rescaling
$$x\to \frac{x}{\sqrt{2}},\;\;\; y\to \frac{y}{\sqrt{2}},$$
and using
$$\frac{dx}{x}=\frac{1}{2}\,\frac{d(x^2)}{x^2},\;\;\;
\theta(x-y)=\theta(x^2-y^2),\;\;{\rm if}\;x>0\;{\rm and}\;y>0,$$
as well as
$$\cos{(x-y)}=[\theta(x-y)+\theta(y-x)]\,\cos{(x-y)},$$
we end up with the result
\begin{equation}
I=K_2=\frac{1}{4\pi^2}\int\limits_0^\infty\frac{dx}{x}\int\limits_0^x
\frac{dy}{y}\,[\cos{(x-y)}-\cos{x}\,].
\label{eq38}
\end{equation}
To calculate the integral
$$\tilde I=\int\limits_0^\infty\frac{dx}{x}\int\limits_0^x
\frac{dy}{y}\,[\cos{(x-y)}-\cos{x}\,],$$
we introduce a parametric integral related to it:
$$\tilde I(\alpha)=\int\limits_0^\infty\frac{dx}{x}\int\limits_0^x
\frac{dy}{y}\,[\cos{(x-\alpha y)}-\cos{x}\,].$$
Note that $\tilde I(0)=0$ and $\tilde I(1)=\tilde I$, so that
$$\tilde I=\int\limits_0^1\frac{d\tilde I(\alpha)}{d\alpha}\,d\alpha.$$
However,
$$\frac{d\tilde I(\alpha)}{d\alpha}=\int\limits_0^\infty\frac{dx}{x}
\int\limits_0^x dy\,\sin(x-\alpha y)=\frac{1}{\alpha}\lim_{\epsilon\to 0}
\int\limits_\epsilon^\infty \frac{dx}{x} \left [\cos{(1-\alpha)x}-\cos{x}\,
\right ],$$
which is the same as
$$\frac{d\tilde I(\alpha)}{d\alpha}=\frac{1}{\alpha}\lim_{\epsilon\to 0}
\left [\int\limits_{\epsilon(1-\alpha)}^\infty \frac{dx}{x}\,\cos{x}
-\int\limits_\epsilon^\infty \frac{dx}{x}\,\cos{x}\right ]=
\lim_{\epsilon\to 0} \frac{Ci(\epsilon)-Ci(\epsilon(1-
\alpha))}{\alpha},$$
where
$$Ci(x)=-\int\limits_x^\infty \frac{dx}{x}\,\cos{x}$$
stands for the integral cosine function. Using the well known series 
representation for this function
$$Ci(x)=\gamma+\ln{x}+\sum\limits_{k=1}^\infty \frac{(-x^2)^k}{2k(2k)!},$$
$\gamma\approx 0.5772$ being the Euler constant, we get
$$\frac{d\tilde I(\alpha)}{d\alpha}=\frac{1}{\alpha}\lim_{\epsilon\to 0}
[\ln{\epsilon}-\ln{(1-\alpha)\epsilon}]=-\frac{\ln{(1-\alpha)}}{\alpha},$$
and, therefore,
$$\tilde I=-\int\limits_0^1\frac{\ln{(1-\alpha)}}{\alpha}\,d\alpha=\zeta(2)
=\frac{\pi^2}{6}.$$
We have just proved an interesting identity which seems to be a new integral
representation for $\zeta(2)$:
\begin{equation}
\zeta(2)=\int\limits_0^\infty\frac{dx}{x}\int\limits_0^x
\frac{dy}{y}\,[\cos{(x-y)}-\cos{x}\,].
\label{eq39}
\end{equation}
Collecting all pieces together, we get finally $I=1/24$ and
$$I_3=\frac{\pi^3}{8}\,\frac{1}{24}=\frac{2}{3!}\left(\frac{\pi}{4}
\right )^3,$$
in agreement with (\ref{eq21}).

\section{Concluding remarks}

This problem had originally aroused in the context of a remarkable
correspondence between the quantum mechanical Landau-Zener problem (known in
the context of molecular scattering) and purely classical problem of a ball
\ (that is 2-sphere) rolling on a Cornu spiral (that is on the curve in 
\textbf{R}$^{2}$, known as Cornu spiral) recently established by Bloch and
Rojo in \cite{2,3}. In fact, the main ingredient of this connection - the
application of the Hopf map - goes back to Feynman, Vernon and Hellwarth 
\cite{4} who showed that the quantum evolution of any two-level system is
determined by the classical evolution (precession) of the magnetic dipole
moment of unit strength in an effective external magnetic field.

It is remarkable that physics helps to calculate a complicated integral
(\ref{eq1}). However, we suspect that there should be a direct method of
calculation. For low values of $n$, we have provided some examples of the 
direct method. It is after the readers to tackle the case of general $n$.

\end{document}